\documentclass[english]{amsart}

\usepackage{babel}
\usepackage{amstext}
\usepackage{amsmath}
\usepackage{amsfonts}
\usepackage{latexsym}
\usepackage{ifthen}
\newcommand\bF{{\mathbb F}}
\newcommand\bR{{\mathbb R}}
\newcommand\bZ{{\mathbb Z}}
\newcommand\bC{{\mathbb C}}
\newcommand\bQ{{\mathbb Q}}

\newcommand\bP{{\mathbb P}}

\newcounter{lemma}

\newtheorem{lemma1}[lemma]{\setcounter{equation}{0}}

\newenvironment{lemma}{\begin{lemma1}{\bf Lemma.}}{\end{lemma1}}

\newenvironment{theorem}{\begin{lemma1}{\bf Theorem.}}{\end{lemma1}}
\newenvironment{question}{\begin{lemma1}{\bf Question.}}{\end{lemma1}}

\newenvironment{proposition}{\begin{lemma1}{\bf Proposition.}}{\end{lemma1}}

\newenvironment{corollary}{\begin{lemma1}{\bf Corollary.}}{\end{lemma1}}

\newenvironment{remark}{\begin{lemma1}{\bf Remark.}\rm}{\end{lemma1}}
\newenvironment{definition}{\begin{lemma1}{\bf Definition.}}{\end{lemma1}}

\newenvironment{Induction Step}{\begin{lemma1}{\bf Induction Step.}}{\end{lemma1}}
\newenvironment{Proof of Theorem 1.2}{\begin{lemma1}{\bf Proof of Theorem 1.2.}}{\end{lemma1}}

\title {The third smallest Salem number in automorphisms of K3 surfaces} 
\author{Keiji Oguiso}
\date{\today}

\begin{document}

\begin{abstract}
We realize the logarithm of the 
third smallest known Salem number as 
the topological entropy of a K3 surface automorphism with a Siegel disk 
and a pointwisely fixed curve at the same time. 
We also show the logarithm 
of the Lehmer number, the smallest known Salem number, 
is not realizable as the topological entropy of 
any Enriques surface automorphism. 
These results are entirely inspired by McMullen's works and Mathematica 
programs.
\end{abstract}

\maketitle

%\tableofcontents

%\vspace*{-0.5in}
\section{Introduction}
\noindent
The aim of this note is to remark the following two new phenomena in complex dymanics of automorphisms of compact complex surfaces. These results and their proofs are entirely inspired by impressive works of McMullen \cite{Mc02-1}, \cite{Mc02-2}, \cite{Mc07}, \cite{GM02} and Mathematica programs. 

\begin{theorem}\label{main} There is a pair $(S, g)$ of a complex K3 surface $S$ and its automorphism $g$ such that:

(1) $S$ contains $8$ smooth rational curves $C_k$ ($0 \le k \le 7$) 
whose dual graph forms the Dynkin diagram $E_8(-1)$ and contains no other 
irreducible complete curve. In particular, $S$ is of algebraic dimension $0$;

(2) The topological entropy $h(g)$ is the logarithm of the third smallest known Salem number 
$$h(g) = \log\, 1.200026523...\, ;$$

(3) The fixed point set $S^g$ consists of one smooth rational curve (in $\cup_{k=0}^{7} C_k$) and $8$ isolated points, say $Q_i$ ($1 \le i \le 7$) and $Q$. The $7$ points $Q_i$ are in $\cup_{k=o}^{7}C_k$, but $Q$ is not in $\cup_{k=0}^{7} C_k$;

(4) $g$ has a Siegel disk at $Q$ and $g$ has no Siegel disk at any other 
point; and 

(5) ${\rm Aut}\, S = \langle g \rangle \simeq \bZ$.   
\end{theorem} 

\begin{theorem}\label{main2} There does not exist a pair $(S, g)$ of a complex Enriques surface $S$ and its automorphism $g$ such that 
$$h(g) = \log\, 1.17628081...\, .$$
Here, the right hand side is the logarithm of the Lehmer number, i.e., the logarithm 
of the smallest known Salem number. 
\end{theorem} 

We shall explain the terms in Theorems (\ref{main}), (\ref{main2}) in Section 2. In the rest of the introduction, we shall remark a few differences between our results and some of preceding known results. 
\par
In \cite{Mc02-2}, McMullen constructed the first examples of surface automorphisms with Siegel disks. They are K3 surface automorphisms arizing from certain Salem numbers of degree $22$, including the $9$-th smallest known one. In his construction, the resultant K3 surfaces are of Picard number $0$. So, they have no complete curve, 
whence, no pointwisely fixed curve as well. Theorem (\ref{main}) tells 
us that it is also possible to have both a Siegel disk and a pointwisely 
fixed curve, necessarily smooth rational, at the same time. 
\par
Let $S$ be a rational surface obtained by blowup at $n$ points on $\bP^2$ 
and $g$ be an automorphism of $S$. Then, $g^{*}(K_S) = K_S$ and $g$ naturally acts on the orthogonal complement 
$K_S^{\perp}$ of the canonical class in $H^2(S, \bZ)$. The lattice $K_S^{\perp}$ is isomorphic to the lattice $E_n(-1)$, i.e., the lattice represented by the Dynkin diagram with $n$ verices $s_k$ ($0 \le k \le n-1$) 
of self-intersection $-2$ such that $n-1$ vertices $s_1$, $s_2$, $\cdots$, $s_{n-1}$ form Dynkin diagram of type $A_{n-1}(-1)$ in this order and the remaining vertex $s_0$ joins to only the vertex $s_3$ by a simple line. (See \cite{Mc07}, Section 2, Figure 2.) The lattice $E_n(-1)$ is of signature $(1, n-1)$ when $n \ge 10$. Then, $g$ naturally induces an orthogonal action $g^{*} \vert E_{n}(-1)$ 
(after fixing a marking). 
By Nagata \cite{Na61} (see also \cite{Mc07}, Theorem (12.4)), 
$g^{*} \vert E_{n}(-1)$ is an element of the Weyl group 
$W(E_{n}(-1))$, i.e., the group generated by the reflections 
$r_k$ ($0 \le k \le n-1$) corresponding to the vertices $s_k$. 
The Weyl group $W(E_n(-1))$ has a special conjugacy class called 
the {\it Coxeter class}. It is the conjugacy class of 
the product (in any order in this case) of the reflections $\Pi_{k=0}^{n-1} r_k$. McMullen (\cite{Mc07}, Theorem (1.1)) shows that, when $n \ge 10$, the Coxeter class is realized {\it algebro-geometrically} by a rational surface automorphism. That is, 
$\Pi_{k=0}^{n-1} r_k = 
g^{*} \vert E_{n}(-1)$ (under a suitable marking) for an automorphism $g$ of $S$ with suitably chosen $n$ blowup points. When $n=10$, i.e., for $E_{10}(-1)$, 
the characteristic polynomial 
of the Coxeter class is exactly the Lehmer polynimial, i.e., the minimal polynomial of the Lehmer number over $\bZ$. In this way, McMullen 
realized the logarithm of the Lehmer number as the topological entropy of some 
rational surface automorphisms with $K_S^{\perp} \simeq E_{10}(-1)$. Note that the Lehmer number is the smallest known Salem number. See \cite{FGR99} and the home page quoted there, for the list of the smallest $47$ known Salem numbers. Being also based on his preceding result \cite{Mc02-1}, Theorem (1.1), 
McMullen (\cite{Mc07}, Theorem (A.1)) also shows that the logarithm of the Lehmer number is in fact the minimal positive entropy of automorphisms of complex surfaces. So, the Lehmer number plays a very special role in automorphisms of compact complex surfaces.
\par
On the other hand, lattice $E_{10}(-1)$ is also isomorphic to the free part 
of $H^2(S, \bZ)$ of an Enriques surface 
$S$. So, it is natural to ask if the logarithm of the Lehmer number can be also realized as the topological entropy of an Enriques surface automorphism or not. Theorem (\ref{main2}) says that it is {\it not}. This may sound negative. However, I believe that such an impossibility result is also of its own interest.
\par
\vskip 4pt
{\it Acknowledgement.} I would like to express my thanks to Professor JongHae Keum for inviting me to the conference ``Algebraic Geometry in East Asia 2008" 
and this volume. This note grew up in my colloquium talk at Osaka University (April 27 2009). I would like to express my 
thanks to Professor Akira Fujiki for inviting me to Osaka University and to Professor Gen Komatsu for giving me an opportunity to talk at the colloquium. This note was completed during my stay at Bayreuth (May 2009 - July 2009) under the renewed program of Alexander von Humboldt Foundation and DFG Forschergruppe 790 ``Classification of algebraic surfaces and compact complex manifolds". Last but not least at all, I would like to express my thanks to Professor Fabrizio Catanese for invitation and staff members here for hospitality and discussions, to Professor De-Qi Zhang for very careful reading. 

%\vfill \eject
%\tableofcontents

\section{Salem numbers in automorphisms of compact K\"ahler surfaces}
\setcounter{lemma}{0}
\noindent
In this section, we quickly review the terms in our Theorems (\ref{main}), (\ref{main2}). As nothing is new, those who are familiar to these terms should 
skip this section.
\par
\vskip 4pt
\noindent
{\it (i) Salem number.} Let us start by the definition.
\begin{definition}\label{defsalem}
A Salem polynomial is a monic irreducible reciprocal polynomial $\varphi(x)$ in $\bZ[x]$ such that $\varphi(x) = 0$ has exactly two real 
roots $\alpha >1$ and $1/\alpha$ off the unit circle $$S^{1} := \{z \in \bC \, \vert\, \vert z \vert \}\,.$$
It is then of even degree. A Salem number is the unique real root $\alpha > 1$. In other words, a Salem number of degree $2n$ is a real algebraic integer $\alpha >1$ whose Galois conjugates 
consist of $1/\alpha$ and $2n-2$ imaginary numbers on $S^{1}$. 
\end{definition}
Salem numbers of degree $2$ are $(m + \sqrt{m^2-4})/2$ ($3 \le m \in \bZ$). For a given integer $n >0$, there are infinitely many Salem numbers of degree $\le 2n$ (\cite{GM02}, Theorem (1.6)). On the other hand, Salem numbers with bounded degree and bounded (Euclidean) norm are finite. That is, for given $n >0$ and $N >0$, Salem numbers $\alpha$ such that ${\rm deg}\, \alpha < 2n$ and $\vert \alpha \vert < N$ are finite. In fact, the elementary symmetric functions of the Galois conjugates of $\alpha$ are then bounded, so that the Salem polynomials of such Salem numbers are finite. So, it is in principle possible to list up all the Salem numbers with explicit boundded norm and 
degree. In fact, there is a list of all Salem numbers of degree $\le 40$ and 
norm $< 1.3$ in the home page quoted by \cite{FGR99}, Page 168. 
The smallest five ones (in degree $\le 40$) are:
$$\alpha_{10} = 1.176280...\, ,\, \alpha_{18} = 1.188368...\, ,$$
$$\alpha_{14} = 1.200026...\, ,\, A_{14} = 1.202616...\, ,$$
$$A_{10} = 1.216391...\, .$$
Their Salem polynomials are
$$\varphi_{10}(x) = x^{10} + x^{9} -x^{7} -x^{6} -x^{5} - x^{4} -x^{3} +x +1\, ,$$$$\varphi_{18}(x) = x^{18} - x^{17} + x^{16} -x^{15} -x^{12} + x^{11} -x^{10} + x^{9} -x^{8} +x^{7} - x^{6} - x^{3} +x^{2} -x  +1\, ,$$ 
$$\varphi_{14}(x) = x^{14} - x^{11} -x^{10} + x^{7} -x^{4} -x^{3} +1\, ,$$
$$\Phi_{14}(x) = x^{14} - x^{12} - x^{7} -x^{2} +1\, ,$$
$$\Phi_{10}(x) = x^{10} -x^{6} -x^{5} - x^{4} +1\, .$$
The smallest Salem number $\alpha_{10}$ (in this range) is called the {\it 
Lehmer number}. This number is discovered by Lehmer \cite{Le33}, Page 477. Lehmer stated there that "{\it We have not made an examination of all $10$th degree symmetric polynomials but a rather intensive search has failed to reveal a better polynomial than $\varphi_{10}(x)$. All efforts to find a better equation of degree $12$ and $14$ have been unsuccessful.}" Since then, it is conjectured that the Lehmer number is the smallest among all the Salem numbers. However, 
it is neither proved nor disproved so far. Also in this view, the result of McMullen \cite{Mc07}, Theorem (A.1) (quoted in the introduction) is very impressive. See also \cite{GH01} and \cite{Mc02-1} for other aspects.
\par
\vskip 4pt
\noindent
{\it (ii) Topological entropy.} Let $X$ be a compact metric space with distance function $d$. Let $g$ be a continuous self map of $X$. To make statement simple, we assume that $g$ is {\it surjective}. The {\it toplogical entropy} is a measure of ``how fast two orbits $\{g^{k}(x)\}_{k \ge 0}$, $\{g^{k}(y)\}_{k \ge 0}$ spread out under $k \rightarrow \infty$". For the precise definition, we introduce a new distance $d_{g, n}$ for each $n \in \bZ_{>0}$ (\cite{KH95}, Page 108):
$$d_{g, n}(x, y) := {\rm max}\, \{d(g^{k}(x), g^{k}(y))\, \vert\, 0 \le k \le n-1\}\, .$$
$d_{g, n}(x, y)$ measures a distance between the orbit segments 
$\{g^{k}(x)\}_{k=0}^{n-1}$, $\{g^{k}(y)\}_{k=0}^{n-1}$. Let $B_{g, n}(x, \epsilon)$ be the open ball 
with center $x$, of radius $\epsilon$ with 
respect to $d_{g, n}$. We call such a ball 
$(\epsilon, n)$-ball. Let 
$S(g, \epsilon, n)$ be the {\it minimal} number of $(\epsilon, n)$-balls 
that cover $X$. Then, ``$S(g, \epsilon, n)$ being larger"
means that the orbit segments of two close points (uniformly with respect to the original distance $d$) 
spread faster in the range $0 \le k \le n-1$. 
\par
\vskip 4pt
The {\it topological entropy} of $g$ is the following value (\cite{KH95}, Page 108, formula (3.1.10)):
$$h(g) := {\rm lim}_{\epsilon \rightarrow 0}\, {\rm limsup}_{n \rightarrow \infty}\, \frac{\log S(g, \epsilon, n)}{n}\, .$$
It is shown that $h(g)$ does not depend on the choice of the distance $d$ giving the same topology on $X$ (\cite{KH95}, Page 109, Proposition 
(3.1.2)). By definition, 
$h(g) = 0$ if $g$ is an automorphism of finite order. 
\par
\vskip 4pt 
Let $E$ be an elliptic curve and $A = E \times E$ be the product abelian surface. By definition, $h(t_{a}) = 0$ for any translation automorphism 
$t_{a}(x) = x + a$ ($a \in A$). Let $M$ be a matrix 
in $M_2(\bZ)$ such that ${\det}\, M \not= 0$. Then $M$ gives rise to the endmorphism $g$ of $A$: $g(x) = Mx$. Let $\alpha$ and $\beta$ be the eigenvalues of $M$ and reorder them so that $\vert \alpha \vert \ge \vert \beta \vert$. Then, according to the three cases
$$(i)\, \vert \alpha \vert \ge \vert \beta \vert \ge 1\,\, , \,\, 
(ii)\, \vert \alpha \vert \ge 1 \ge \vert \beta \vert\,\, ,\,\,  
(iii)\, 1 \ge \vert \alpha \vert \ge \vert \beta \vert\, ,$$
we have
$$(i)\, h(g) = \log\, \vert \alpha \beta \vert^2\, ,\, (ii)\, 
h(g) = \log\, \vert \alpha \vert^2\,\, ,\,\, (iii)\, h(g) = \log 1 = 0\, .$$
The reason is as follows. First, choose sufficiently small $\epsilon >0$ and cover $A$ by $N$ mutually disjoint complex $2$-dimensional $\epsilon$-cubes (in the product metric) that are "parallel to" the two complex eigenvectors of $M$ (
Here we ignore small part of $A$ in covering). Next, divide each of $N$ $\epsilon$-cubes into mutually disjoint $(\epsilon, n)$-cubes with respect to the new distance $d_{g, n}$. Then, according to the cases (i), (ii), (iii), the numbers of the resultant $(\epsilon, n)$-cubes are approximately 
$N \cdot \vert \alpha \beta \vert^{2(n-1)}$, 
$N \cdot \vert \alpha\vert^{2(n-1)}$, 
and $N$. This implies the result. (See \cite{KH95}, Pages 121-123, for more 
precise calculations). The values $h(g)$ above coincide with the logarithm of the spectral radius of 
the action of $g^{*} \vert \oplus_{k=0}^{2} H^{2k}(A, \bZ)$. However, 
this is not accidental:
\begin{theorem}\label{gyf} 
Let $X$ be a compact K\"ahler manifold of dimension $n$ 
and let $g : X \longrightarrow X$ 
be a holomorphic surjective self map of $X$. Then 
$$h(g) = \log \rho(g^{*} \vert \oplus_{k=0}^{n} H^{2k}(X, \bZ))\, .$$
Here $\rho(g^{*} \vert \oplus_{k=0}^{n} H^{2k}(X, \bZ))$ is 
the spectral radius 
of the action of $g^{*}$ on the total cohomology ring of even degree.
\end{theorem}
This is a fundamental theorem often attributed to Gromov and Yomdin. 
The explicit statement with full proof (using Yomdin's result) 
is found in Friedland's paper \cite{Fr95}, Theorem (2.1). Note that, in the 
proof, we only need the estimate by the spectral radius on the cohomology group of even degree. See also \cite{DS04}, Pages 315-316, for further discussions.
As an immediate consequence, we obtain the following important
\begin{corollary}\label{alg} 

(1) $h(g)$ is the logarithm of an algebraic integer. 

(2) $h(g^n) = nh(g)$ for a positive integer $n$.
\end{corollary}
\par
\vskip 4pt
\noindent
{\it (iii) Toplogogical entropy of a surface automorphism.} If $\dim\, X = 1$, 
then by Theorem (\ref{gyf}), $h(g) = \log\,({\rm deg}\, g)$ and it is not so informative. Let us consider the case where $X$ is a compact K\"ahler surface 
and $g$ is an {\it automorphism} of $X$. 
\par
\vskip 4pt
The first important fact 
is the following result due to Cantat (\cite{Ca99}, Proposition 1 and 
its proof):
\begin{theorem}\label{ca}
Let $X$ be a compact K\"ahler surface and $g$ be an automorphism of $X$. 
Then, 

(1) $h(g) = \log\, \rho(g^{*} \vert H^{1,1}(S, \bR))$. Here $\rho(g^{*} \vert H^{1,1}(S, \bR))$ is the spectral radius of $g^{*} \vert H^{1,1}(S, \bR)$.

(2) Assume that $h(g) > 0$. Then $X$ is isomorphic to either:

(i) a rational surface with $b_2(X) \ge 11$; 

(ii) a $2$-dimensional complex torus (or its blow up)

(iii) a K3 surface (or its blow up); or

(iv) an Enriques surface (or its blow up).
\end{theorem}
See eg. \cite{BHPV04}, Pages 244-246, for the classification of compact complex surfaces and the definition of surfaces above. Non-minimality in (i) 
is definitely essential. But non-minimal surfaces in (ii)-(iv) are not 
so essential. 
In fact, if $X$ is not minimal in the class (ii)-(iv), then 
$g$ descends to the automorphism $\overline{g}$ of the minimal model 
(See eg. \cite{BHPV04}, Page 99, Claim). Moreover, 
$h(g) = h(\overline{g})$. This follows from 
the fact that the exceptional set forms a negative definite sublattice of 
$H^{2}(X, \bZ)$. So, surfaces having interesting automorphisms 
(in the view of topological entropy) 
are non-miminal rational surfaces, $2$-dimensional complex tori, 
K3 surfaces and Enriques surfaces.
\par
\vskip 4pt
Salem numbers naturally appear in the study of automorphisms of 
complex surfaces:
\begin{theorem}\label{factorization}
Let $X$ be a compact K\"ahler surface and $g$ be an automorphism of $X$. 
Then the characteristic polynomial of $g^{*} \vert H^2(X, \bZ)$ is the product 
of cyclotomic polynomials and Salem polynomials. In the product, 
there are at most one Salem factor (counted with multiplicities) and 
possibly no cyclotomic factor or no Salem factor. In particular, if $h(g) >0$, then Salem factor appears in the product and $h(g)$ is the logarithm of that Salem number. 
\end{theorem}
This is due to McMullen \cite{Mc02-2}, Theorem (3.2). The argument 
there is given for K3 surface automorphisms. But it is easily generalized to 
automorphisms of arbitrary compact K\"ahler surfaces. 
\begin{proof} Consider the real Hodge decomposition of $H^2(X, \bZ)$:
$$H^{2}(X, \bR) = H_{\bR}^{1,1}(X) \oplus V\, .$$
Here $V$ is a vector subspace of $H^{2}(X, \bR)$ such that
$$V \otimes \bC = H^{2,0}(X) \oplus H^{0,2}(X)\, .$$
The K\"ahler cone ${\mathcal K}(X)$ forms a strictly convex open cone of 
$H_{\bR}^{1,1}(X)$. Moreover, $H_{\bR}^{1,1}(X)$ is of signature $(1, h^{1,1}(X) -1)$ 
and $V$ is positive definite (see eg. \cite{BHPV04}, Page 143, Theorem (2.14)).
As $g$ is an automorphism, we have $g^{*}(V) = V$, $g^{*}(H_{\bR}^{1,1}(X)) 
= H_{\bR}^{1,1}(X)$ 
and $g^{*}({\mathcal K}(X)) = 
{\mathcal K}(X)$. As $V$ is negative definite, 
the eigenvalues of $g^{*} \vert V$ are of absolute value $1$. 
As $g^{*} \in {\rm O}(H^2(X, \bR))$, it follows that 
${\rm det}\, g^{*} \vert H^{2}(X, \bR) = \pm 1$. 
Thus, the product of the eigenvalues of $g^{*} \vert H_{\bR}^{1,1}(X)$ 
is of absolute value $1$ as well. On the other hand, as 
$g^{*} \vert H_{\bR}^{1,1}(X)$ preserves the strictly convex open cone 
${\mathcal K}(X)$, the spectral radius of $g^{*} \vert H_{\bR}^{1,1}(X)$ 
is given by a real eigenvalue, say $\alpha >0$,  
of $g^{*} \vert H_{\bR}^{1,1}(X)$ 
with eigenvector, say $\eta$, in $\overline{{\mathcal K}(X)}$ (the closure 
of the K\"ahler cone). 
This is due to the (generalized) 
Perron-Frobenious theorem (\cite{Bi67}, Page 274, Theorem). 
As the product of the eigenvalues is of absolute value $1$, it follows 
that $\alpha \ge 1$. If $\alpha = 1$, then all the eigenvalues of $g^{*} \vert H_{\bR}^{1,1}(X)$ is also of absolute value $1$. Hence, so are the eigenvalues 
of $g^{*} \vert H^{2}(X, \bR)$. As $g^{*} \vert H^{2}(X, \bR)$ is defined over 
$H^{2}(X, \bZ)$, all the eigenvalues are then roots of unity 
by Kronecker's theorem. Next consider the case where $\alpha >1$. Consider 
$g^{-1}$. Then $1/\alpha <1$ is an eigenvalue of $(g^{-1})^{*} \vert H_{\bR}^{1,1}(X)$. Then, again, by the generalized Perron-Frobenious, the spectral radius of $(g^{-1})^{*} \vert H_{\bR}^{1,1}(X)$ is a real eigenvalue, say $\beta$, 
with eigenvector, say $\eta'$, in $\overline{{\mathcal K}(X)}$. 
By $1/\alpha <1$, 
we have $\beta >1$. Thus, $g^{*} \vert H_{\bR}^{1,1}(X)$ 
has an eigenvalue $\alpha' := 1/\beta$ with the same eigenvector $\eta'$. 
As $\alpha \not= \alpha'$, the linear subspace $H = \bR\langle \eta, \eta' \rangle$ is $2$-dmensional. From 
$$(g^{*}\eta, g^{*}\eta) = \alpha^2(\eta, \eta)\,\, ,\,\, (g^{*}\eta', 
g^{*}\eta') = (\alpha')^2(\eta', \eta')\, ,$$
we obtain $(\eta, \eta) = (\eta', \eta') = 0$. Moreover, we have 
$(\eta, \eta') > 0$, as both vectors are in the closure of the K\"ahler cone 
$\overline{{\mathcal K}(X)}$ (and are linearly independent). Thus, $H$ is of signature $(1,1)$ and the orthogonal complement $V^{\perp}$ in $H_{\bR}^{1,1}(X)$ 
is negative definite. Thus, the remaining eigenvalues of $g^{*} \vert 
H_{\bR}^{1,1}(X)$, that coincide with the eigenvalues of $g^{*} \vert V^{\perp}$, are of absolute value $1$. In conclusion, 
$g^{*} \vert H^{2}(X, \bR)$ has two real 
eigenvalues $\alpha >1$, $0 < \alpha' <1$ and the other eigenvalues are all 
of absolute value $1$. Also, all these values are algebraic integers. 
This is because $g^{*}\vert H^{2}(X, \bR)$ is defined over $H^{2}(X, \bZ)$. This implies 
the result.   
\end{proof}
Finally we recall the notion of a Siegel disk (for simplicity only 
$2$-dimensional case). 
\begin{definition}\label{sd}

(1) Let $\Delta^2$ be a $2$-dimensional unit disk 
with linear coordinate $(z_1, z_2)$. A linear automorphism (written under the coordinate action)
$$f^{*}\left( \begin{array}{c} z_1\\ 
z_2 \end{array}\right)\ = \left( \begin{array}{cc} \rho_1 & 0\\ 
0 & \rho_2 \end{array}\right)\ \left( \begin{array}{c} z_1\\ 
z_2 \end{array}\right)\ \,$$
is called an irrational rotation if $\vert \rho_1 \vert = \vert \rho_2 \vert = 1$, and $\rho_1$ and $\rho_2$ are multiplicatively independent, in the sense 
that 
$(m_1, m_2) = (0, 0)$ is the only integer solution to 
$$\rho_1^{m_1}\rho_2^{m_2} = 1\, .$$

(2) Let $X$ be a complex analytic surface (not necessarily compact) 
and $g$ be an automorphism of $X$. A domain $U \subset X$ is called 
a Siegel disk of $(S, g)$ if $g(U) =U$ and $(U, g \vert U)$ 
is isomorphic to some irrational rotation $(\Delta^2, f)$. In other words, 
$g$ has a Siegel disk if and only if there is a fixed point $P$ at which 
$g$ is locally 
analytically linearlized as in the form of an irrational rotation. 
\end{definition}
The existence of a Siegel disk implies that 
there is no topologically dense orbit. The first examples of 
surface automorphisms with Siegel 
disks were discovered by McMullen (\cite{Mc02-2}, Theorem (1.1)) within 
K3 surfaces. 
The resultant K3 surfaces $X$ 
are necessarily of 
algebraic dimension $0$. This follows from the fact that 
the action on the space of holomorphic $2$-forms is finite cyclic if 
the algebraic dimension $\not= 0$ (\cite{Mc02-2}, Theorem (3.5), see also 
\cite{Og08}, Theorem (2.4)). We also note that in this case $NS(X)$ 
is negative definite (and vice versa), so that $X$ contains at most finitely 
many irreducible complete curves 
and they are all smooth rational (if exist). This easily follows from 
the Riemann-Roch 
inequality for K3 surfaces (see eg. \cite{BHPV04}, Page 312, line 6, 
formula (2)). Later, 
McMullen (\cite{Mc07}, Theorem (10.1)) also found rational surface automorphisms with Siegel disks. In this case, 
the resultant surfaces are projective. In fact, they are blowup of 
$\bP^2$.
\par
\vskip 4pt 
In general, it is hard to see if a given action is locally analytically linearizable at the fixed point or not. The following criterion, which we only state in dimension $2$, is again due to 
McMullen (\cite{Mc02-2}, Theorem (5.1)): 
\begin{theorem}\label{siegeldisk} 
Let $\varphi$ be an automorphism of a germ of the origin $0$ of $\bC^2$ 
such that $\varphi(0) = 0$ and such that 
$$d\varphi^{*}(0) = \left( \begin{array}{cc} \rho_1 & 0\\ 
0 & \rho_2 \end{array}\right)\ \, .$$
Here $d\varphi^{*}(0)$ is the action on the cotangent space 
$\Omega_{\bC^2}^{1}(0)$ at $0$ induced by the coordinate action $\varphi^{*}$. 
(We prefer coordiante action as then everything is covariant.) Assume that: 

(1) $\rho_1$ and $\rho_2$ are algebraic numbers;

(2) $\vert \rho_1 \vert = \vert \rho_2 \vert = 1$; and

(3) $\rho_1$ and $\rho_2$ are multiplicatively independent. 

Then $\varphi$ has a Siegel disk at $0$, i.e., there is a local coordinate 
$(z_1, z_2)$ at $0$ such that 
$$\varphi^{*}\left( \begin{array}{c} z_1\\ 
z_2 \end{array}\right)\ = \left( \begin{array}{cc} \rho_1 & 0\\ 
0 & \rho_2 \end{array}\right)\ \left( \begin{array}{c} z_1\\ 
z_2 \end{array}\right)\ \, .$$
\end{theorem}
This is a highly non-trivial result that involves very deep theorems: 
the Siegel-Sternberg theorem on 
analytic linearlization and the Baker-Fel'dman theorem on transcendence of 
the logarithm of algebraic numbers. See \cite{Mc02-2}, Section 5 and the references therein for more details. 
\par
\vskip 4pt

\section{Proof of Theorem (\ref{main})}
\setcounter{lemma}{0}
\noindent
Let us consider the Salem polynomial
$$\varphi_{14}(x) = x^{14} - x^{11} - x^{10} +x^{7} - x^{4} -x^{3} +1\, .$$
The third smallest Salem number 
$$\alpha_{14} = 1.200026...$$
is the unique real root $> 1$ of $\varphi_{14}(x) =0$. 
The equation $\varphi_{14}(x) =0$ 
has one more real root $1/\alpha_{14}$. The other $12$ roots, which we 
denote by 
$$\beta_k\,\, ,\,\, \overline{\beta_k}\,\, (1 \le k \le 6)\, ,$$
are on the unit circle $S^1 = \{z \in \bC\, \vert\, \vert z \vert 
= 1\, \}$. Among these $12$ roots on $S^1$, we choose two particular ones: 
$$\delta := \beta_1 := -(0.990398835230041...) - (0.31823945592693...)i\, ,$$
$$\theta := \beta_2 := -(0.371932997164175...) - (0.92825957879273...)i\, .$$
These approximate values are computed by Mathematica program, NSolve.
\par
\vskip 4pt
{\it In what follows, $\delta$ and $\theta$ always mean these 
two particular roots.}
\par
\vskip 4pt 
We denote the K3 lattice by:
$$\Lambda := \Lambda_{\rm K3} := E_{8}(-1)^{\oplus 2} \oplus 
H^{\oplus 3}\, .$$
Here $H$ is the unique even unimodular lattice of signature 
$(1,1)$ and $E_8(-1)$ is the unique even unimodular negative definite lattice 
of rank $8$. The lattice $\Lambda$ is isomorphic to the second cohomology lattice $(H^2(S, \bZ), (*,**))$ of a K3 surface $S$. Here $(*,**)$ is the cup product on $H^2(S, \bZ)$. (See eg. \cite{BHPV04}, Page 311, Proposition (3.3)(ii).)  For a field $K$, we denote the $K$-vector space $\Lambda \otimes_{\bZ} K$ 
by $\Lambda_{K}$. A similar abbreviation will be applied for other lattices and vector spaces.
\begin{proposition}\label{lattice}
There are an automorphism $F$ of the K3 lattice $\Lambda$
and an element $\sigma$ of $\Lambda_{\bC}$ such that: 

(1) $\Phi_F(x) = (x-1)^8 \cdot \varphi_{14}(x)$, where $\Phi_F(x)$ 
is the characteristic polynomial of $F$;

(2) $(\sigma , \sigma) = 0$ and $(\sigma, \overline{\sigma}) > 0$; and

(3) $F(\sigma) = \delta \sigma$.
\end{proposition}

\begin{proof} Recall the following theorem \cite{GM02}, Theorem 1.3.

\begin{theorem}\label{GM} Let $\varphi(x) \in \bZ[x]$ be an irreducible reciprocal polynomial such that $\vert \varphi(\pm 1) \vert = 1$ and $p$, $q$ 
be positive integers such that $p \equiv q\, ({\rm mod}\, 8)$. Let $\bR^{p+q}$ be the real vector space with a symmetric bilinear from of signature 
$(p, q)$. Assume that $f \in {\rm SO}(\bR^{p+q})$ and 
the characteristic polynomial $\Phi_{f}(x)$ is $\varphi(x)$. Then, there is an even unimodular lattice $L \subset \bR^{p+q}$ such that $L_{\bR} = \bR^{p+q}$ and $f(L) = L$. In other words, $f$ is realized as an automorphism of an even unimodular lattice of signature $(p,q)$.  
\end{theorem} 
Our $\varphi_{14}(x)$ is an irreducible reciprocal polynomial in $\bZ[x]$ with $\vert \varphi_{14}(\pm 1) \vert = 1$. We apply Theorem (\ref{GM}) 
for $\varphi_{14}(x)$, along the line of \cite{GM02}, Pages 270--271, Proof 
of Theorem 2.2. 
\par
\vskip 4pt
Let $V_{k}$ ($0 \le k \le 6$) be the 
real vector space 
$\bR^2$. We define symmetric bilinear forms $Q_k$ on $V_k$ ($0 \le k \le 6$) 
by 
$$Q_0 := \left( \begin{array}{cc} 0 & 1\\ 
1 & 0 \end{array}\right)\ \, ,\, 
Q_1 := I_2 :=  \left( \begin{array}{cc} 1 & 0\\ 
0 & 1 \end{array}\right)\ \, ,\,
Q_k := -I_2 :=  \left( \begin{array}{cc} -1 & 0\\ 
0 & -1 \end{array}\right)\ \, (k \ge 2)\, .$$
Then $(V_0, Q_0)$ is of signature $(1,1)$, $(V_1, Q_1)$ is positive definite 
and $(V_k, Q_k)$ ($k \ge 2$) are negative definite. 
Noticing $\vert \beta_k \vert = 1$, we define $f_k \in {\rm SO}(V_k, Q_k)$ ($0 \le k \le 6$) by 
$$f_0 = \left( \begin{array}{cc} \alpha_{14} & 0\\ 
0 & 1/\alpha_{14} \end{array}\right)\ \, ,\, 
f_k := \left( \begin{array}{cc} {\rm Re}\, \beta_k & -{\rm Im}\, \beta_k\\ 
{\rm Im}\, \beta_k & {\rm Re}\, \beta_k \end{array}\right)\  \, (k \ge 1)\, .$$ Here $\alpha_{14}$, $1/\alpha_{14}$, $\beta_k$, $\overline{\beta_k}$ ($1 \le k \le 6$) are the roots of $\varphi_{14}(x) = 0$. The eigenvalues of 
$f_0$ are $\alpha_{14}$ and $1/\alpha_{14}$, and the eigenvalues of $f_{k}$ 
($k \ge 1$) are $\beta_k$ and $\overline{\beta_k}$. Set 
$$(V, Q) := (\oplus_{i=0}^{6}V_i, \oplus_{i=0}^{6}Q_i)\, ,\, f := 
\oplus_{k=0}^{6} f_k\, .$$ 
By construction, $(V, Q) = \bR^{3+11}$, $f \in {\rm SO}(\bR^{3+11})$ and 
the characteristic polynomial of $f$ is $\varphi_{14}(x)$.  
Thus, by Theorem (\ref{GM}), there is an even unimodular lattice 
$L \subset \bR^{3+11}$ such that $f(L) = L$ and $L_{\bR} = \bR^{3+11}$. 
We have an isomorphism
$$L \simeq E_8(-1) \oplus H^{\oplus 3}\, .$$
This is because the isomorphism class of 
an even indefinite unimodular lattice is uniquely determined by 
the signature (\cite{Se73}, Page 54, Theorem 5). 
\par
\vskip 4pt 
We can thus identify  
$$\Lambda = E_{8}(-1) \oplus L\, .$$
Put $F = id_{E_{8}(-1)} \oplus f$. Then $F \in {\rm SO}(\Lambda)$ 
and the characteristic polynomial of $F$ is 
$$\Phi_{F}(x) = (x-1)^8 \cdot \varphi_{14}(x)\, .$$ 

It remains to find $\sigma \in \Lambda_{\bC}$ that satisfies 
(2) and (3). Choose 
an eigenvector $\sigma \in \Lambda_{\bC}$ of $F$ with eigenvalue 
$\delta = \beta_1$. We shall show that this $\sigma$ satisfies 
(2) and (3). By definition, we have $F(\sigma) = \delta \sigma$. As 
$F$ is an automorphism of the {\it lattice} $\Lambda$, 
it follows that  
$$(\sigma, \sigma) = (F(\sigma), F(\sigma)) = \delta^2 (\sigma, \sigma)\, .$$
Thus, $(\sigma, \sigma) = 0$ by $\delta^2 \not= 1$. Taking the complex 
conjugate, we obtain that 
$$F(\overline{\sigma}) = \overline{\delta}\overline{\sigma}\, ,\, 
(\overline{\sigma}, \overline{\sigma}) = 0\, .$$
Note that  
$$\sigma + \overline{\sigma} \not= 0\, .$$ 
This is because $\sigma$ and $\overline{\sigma}$ are eigenvectors with different eigenvalues. On the other hand, by the explicit form of $F$, we see that 
$$\sigma\, ,\, \overline{\sigma} \in  (V_1)_{\bC} \subset L_{\bC} \subset 
\Lambda_{\bC}\, .$$  
As $Q_1$ is positive definite on $V_1$ and $\sigma + \overline{\sigma}$ is 
a real vector in $V_1 \setminus \{0\}$, it follows that 
$$(\sigma + \overline{\sigma}, \sigma + \overline{\sigma}) > 0\, .$$ 
As $(\sigma, \sigma) =(\overline{\sigma}, \overline{\sigma}) = 0$, 
this implies $(\sigma, \overline{\sigma}) >0$. 
\end{proof}

\begin{remark}\label{lattice2} By changing the symmetric bilinear 
forms 
on $V_1$ by $-I_2$ and on $V_2$ by $I_2$, we have an automorphism $F'$ 
of the K3 lattice $\Lambda$ and an element $\sigma' \in \Lambda_{\bC}$ 
such that 

(1) $\Phi_{F'}(x) = (x-1)^8 \cdot \varphi_{14}(x)$;

(2) $(\sigma' , \sigma') = 0$ and $(\sigma', \overline{\sigma'}) > 0$; and

(3) $F'(\sigma') = \theta \sigma'$ (Here we recall that $\theta = \beta_2$).
\end{remark}

\begin{theorem}\label{K3} There is a pair $(S, g)$ 
of a K3 surface $S$ and its automorphism $g$ 
such that 

(1) $g^*\sigma_S = \delta \sigma_S$; 

(2) The N\'eron-Severi lattice $NS(S)$ is isomorphic to $E_8(-1)$; and 

(3) $g^* \vert NS(S) = id_{NS(S)}$. 
\end{theorem}

See eg. \cite{BHPV04}, Page 308, line 4 for the definition of the N\'eron-Severi lattice $NS(S)$. 

\begin{proof} Let $F$ and $\sigma$ be the same as in 
Proposition (\ref{lattice}). 
Then, by Proposition (\ref{lattice}) (2), the point $[\bC \sigma]$ belongs 
to the period domain of K3 surfaces:
$$\Omega := \{[\bC \sigma] \in \bP(\Lambda_{\bC}) \, \vert\, (\sigma, \sigma) = 0\, ,\, (\sigma, 
\overline{\sigma}) > 0\}\, .$$ Thus, we can apply 
the surjectivity of the period 
mapping for K3 surfaces (see eg. \cite{BHPV04}, Page 339, Corollary (14.2)) 
to get a K3 surface $S$ 
and an isomorphism $\iota : H^2(S, \bZ) \simeq \Lambda$ such that 
$\iota(\bC \sigma_S) = \bC \sigma$.  
Here $H^0(S, \Omega_S^2) = \bC \sigma_S$. Define an automorphism $f_S$ of 
$H^2(S, \bZ)$ by  
$$f_S := \iota^{-1} \circ F \circ \iota\, .$$ 
We want to find an automorphism $g$ of $S$ such that 
$f_S = g^{*}$. According to the global Torelli theorem for K3 surfaces (see eg. \cite{BHPV04}, Page 332, Theorem 
(11.1)), this follows if $f_S$ satisfies the following three properties 
(i)-(iii): 

(i) $f_S$ is an Hodge isometry;

(ii) $f_S$ preserves the positive cone $\mathcal P(S)$, i.e., the connected component of 
$$\{x \in H^{1,1}(S, \bR)\, \vert (x, x) > 0\}$$
containing the K\"ahler classes of $S$; and

(iii) $f_S$ preserves the set of classes represented by effective curves in $NS(S)$. 
\par
\vskip 4pt
Let us check these properties. By definition of $f_S$, we have 
$f_S \in {\rm SO}\,(H^2(S, \bZ))$ and 
$f_S(\sigma_S) = \delta \sigma_S$. This shows (i). Recall that the Salem number $\alpha_{14}$ is real and an eigenvalue of $F$. So, $\alpha_{14}$ 
is a real eigenvalue of $f_{S}$ as well. We can then choose a real eigenvector
 $\eta \in H^2(S, \bR)$ of $f_S$ with eigenvalue $\alpha_{14}$. By 
$$(\eta, \sigma_S) = (f_S^{*}(\eta),f_S^{*}(\sigma_S)) = \alpha_{14}\delta 
(\eta, \sigma_S)$$ 
and by $\alpha_{14}\delta \not= 1$, 
we have $(\eta, \sigma_{S}) = 0$. As $\eta$ is real, this implies that 
$\eta \in H^{1,1}(S, \bR)$. Moreover, by $\alpha_{14} >1$ and by 
$$(\eta, \eta) = (f_S^{*}(\eta),f_S^{*}(\eta)) = 
\alpha_{14}^2(\eta, \eta)\, ,$$ 
we have $(\eta, \eta) = 0$. Thus, $\eta \in \partial {\mathcal P}(S)$ (the boundary of the positive cone) possibly 
after replacing $\eta$ by $-\eta$. 
As $f_S(\eta) = \alpha_{14}\eta$ with $\eta \not= 0$ and $\alpha_{14} > 0$, 
this implies (ii). 
\par
\vskip 4pt
It remains to check (iii). The $\bC$-linear extension of the lattice $L \simeq 
E_8(-1) \oplus H^{\oplus 3}$ (defined in the proof of Proposition 
(\ref{lattice})) contains $\sigma$. Moreover, as the characteristic polynomial 
of $F \vert L = f$, which is $\varphi_{14}(x)$, is irreducible over 
$\bZ$, it follows that the lattice $L$ is 
the {\it minimal} primitive lattice of 
$\Lambda$ of which $\bC$-linear extension contains $\sigma$. Thus, 
the lattice $\iota^{-1}(L)$ is also the minimal primitive sublattice of $H^2(S, \bZ)$ of which $\bC$-linear extension contains 
$\sigma_S$. By definition of the transcendental lattice $T(S)$ (See eg. 
\cite{BHPV04}, Page 308, line 5), we have then 
that  
$$T(S) = \iota^{-1}(L) \simeq E_8(-1) \oplus H^{\oplus 3}\, .$$ 
Recall that $NS(S) = T(S)^{\perp}$ in $H^{2}(S, \bZ)$ and 
$L^{\perp} = E_{8}(-1)$ in $\Lambda$. Then 
$$NS(S) = \iota^{-1}(L^{\perp}) = \iota^{-1}(E_8(-1)) 
\simeq E_8(-1)\, .$$
Moreover, by $F = id_{E_8(-1)} \oplus f$, it follows that 
$$f_S(T(S)) = T(S)\, ,\, f_S(NS(S)) = NS(S)\,\, {\rm and}\,\, 
f_S \vert NS(S) = id_{NS(S)}\, .$$ 
So, the assertion (iii) holds. 
\par
\vskip 4pt
Hence there is an automorphism $g$ of $S$ such that 
$f_S = g^{*}$. By the proof of (iii), 
our $(S, g)$ also satisfies the assertions (2) and (3) of Theorem (\ref{K3}). 
This completes the proof. 
\end{proof}

\begin{remark}\label{K3'} Starting from $F'$ and $\sigma'$ in Remark (\ref{lattice2}) 
(instead of $F$ and $\sigma$ in Proposition (\ref{lattice})), we also 
obtain a pair 
$(S', g')$ of a K3 surface $S'$ and its automorphism $g'$ such that 

(1) $(g')^*\sigma_{S'} = \theta \sigma_{S'}$; 

(2) $NS(S') \simeq E_8(-1)$; and 

(3) $(g')^* \vert NS(S') = id_{NS(S')}\,$. 
\end{remark}

{\it In the rest, we shall show that the pair $(S, g)$ satisfies 
the requirement 
of Theorem (\ref{main}) but the pair $(S', g')$ does not.}

\begin{proposition}\label{ns} Let $S$ be a K3 surface such that 
$NS(S) \simeq E_8(-1)$. Then, $S$ contains $8$ smooth rational curves 
$C_k$ ($0 \le k \le 7$) and contains no other irreducible complete 
curve. Moreover, the dual graph 
of $C_k$ ($0 \le k \le 7$) is the same as the Dynkin diagram $E_8(-1)$, 
i.e., $(C_k^2) = -2$, vertices $C_1$, $C_2$, $\cdots$, $C_7$ form Dynkin diagram of Type $A_7(-1)$ in this order and the vertex $C_0$ joins to only the vertex $C_3$ by a simple line. (See \cite{Mc07}, Section 2, Figure 2. In the figure, set $n=7$ and replace $s_k$ there by $C_k$ here.)
\end{proposition}

\begin{proof} We shall show Proposition (\ref{ns}) by dividing into four steps.
\par
\vskip 4pt
\noindent
{\it Step 1.} {\it Let $C$ be an irreducible complete curve on $S$. Then 
$C \simeq \bP^1$.}
\par
\vskip 4pt
\noindent
{\it Proof.} As $NS(S)$ is even, negative definite and $S$ is K\"ahler (see eg. \cite{BHPV04}, Page 144, Theorem (3.1) and Page 310 Proposition (3.3)(i)), 
we have that $C \not\equiv 0$ and $(C^2) \le -2$. Thus, for 
the arithmetic genus $p_a(C)$, we have 
$$0 \le p_a(C) = (C^2)/2 + 1 \le 0\, .$$
Hence $p_a(C) = 0$. This implies $C \simeq \bP^1$. \qed
\par
\vskip 4pt
\noindent
{\it Step 2.} {\it $NS(S)$ is generated by the classes of irreducible 
complete curves. In particular, the number of irreducible complete 
curves on $S$ 
is greater than or equal to $8$.}
\par
\vskip 4pt
\noindent
{\it Proof.} Let $e_i$ ($0 \le i \le 7$) be the basis of $NS(S)$ corresponding 
to the $8$ vertices of $E_8(-1)$. We have $(e_i^{2}) = -2$. Let $E_i \in {\rm Pic}\, S$ be a representative of $e_i$. Then by the Riemann-Roch formula and the Serre duality, we have 
$$h^0(E_i) + h^0(-E_i) \ge \frac{(E_i^2)}{2} + 2 = 1\, .$$
Thus, for each $i$, either $\vert E_i \vert$ or $\vert -E_i \vert$ contains 
an effective curve. As the class of each irreducible component is also in $NS(S)$, this implies the result. \qed
\par
\vskip 4pt
\noindent
{\it Step 3.} {\it Let $C_k$ ($0 \le k \le m-1$) be mutually distinct irreducible complete curves on $S$. Then the classes $[C_k] \in NS(S)$ are linearly independent 
in $NS(S)$. In particular, the number of irreducible complete curves on $S$ is less than or equal to $8$.}
\par
\vskip 4pt
\noindent
{\it Proof.} If otherwise, there are subsets $I$ and $J$ of $\{0, 1, \cdots , m-1\}$ such that $I \cap J = \emptyset$ and 
$$\sum_{i \in I} a_i[C_i] = \sum_{j \in J} b_j[C_j]\, .$$
Here $a_i \ge 0$ and $b_j \ge 0$ and $a_i \not= 0$ for at least one $a_i$. 
As $NS(S)$ is negative definite, it follows that 
$$0 > ((\sum_{i \in I} a_i[C_i])^2)\, .$$
On the other hand, we have that
$$((\sum_{i \in I} a_i[C_i])^2) = 
(\sum_{i \in I} a_i[C_i], \sum_{j \in J} b_j[C_j] ) \ge 0\, ,$$
a contradiction. This implies the result. \qed
\par
\vskip 4pt
\noindent
{\it Step 4.} {\it $S$ contains $8$ smooth rational curves 
whose dual graph forms Dynkin diagram $E_8(-1)$ and contains no other irreducible complete 
curve.}
\par
\vskip 4pt
\noindent
{\it Proof.} By Steps 2,3, $S$ contains exactly $8$ irreducible complete 
curves. We denote them by  
$C_k$ ($0 \le k \le 7$). Again by Steps 2,3, 
$\langle [C_k] \rangle_{k=0}^{7}$ form a basis of $NS(S)$ 
over $\bZ$. By Step 1, each $C_k$ is also a smooth rational curve. Thus 
$(C_k^2) = -2$. As $NS(S)$ is negative definite, the dual graph of $\{C_k\}_{k=0}^{7}$ is then a disjoint union of 
Dynkin diagrams of type $A_n(-1)$, $D_m(-1)$, $E_6(-1)$, $E_7(-1)$, $E_8(-1)$, 
with $8$ vertices in total. As $NS(S)$ is unimodular and $\langle [C_k] \rangle_{k=0}^{7}$ 
forms a basis of $NS(S)$ over $\bZ$, the only possible dual graph 
of $\{C_k\}_{k=0}^{7}$ is then $E_8(-1)$. 
In fact, 
the lattices associated with other Dynkin diagrams are 
of discriminant $\ge 2$. This completes the proof. 
\end{proof}

Let us return back to our $(S, g)$ in Theorem 
(\ref{K3}). $S$ has exactly $8$ smooth rational curves, say $C_k$ ($0 \le k \le 7$), as described in Proposition (\ref{ns}), and no other irreducible complete curve. We set $S^g := \{x \in S\, \vert\, g(x) = x\}$.

\begin{lemma}\label{fixcurve} 

(1) $g(C_k) = C_k$ for each $C_k$.

(2) Put $\mathcal F_C = S^g \cap (\cup_{k=0}^{7} C_k)$. Then, 
$$\mathcal F_C = C_3 \cup \{P_1, P_{12}, P_0, P_{45}, P_{56}, 
P_{67}, P_7\}\, .$$
Here $P_{ij}$ is the intersection point of $C_i$ and $C_{j}$, and $P_i$ 
is a point on $C_i \setminus \cup_{j \not= i} C_j$. Moreover, for each $P \in 
{\mathcal F}_C$, the action $dg^{*}(P)$ of $g^{*}$ (the coordinate action of 
$g$) on the cotangent space 
$\Omega_{S}^1(P)$ is diagonalized as follows:

$$dg^{*}(P) = \left( \begin{array}{cc} 1 & 0\\ 
0 & \delta \end{array}\right)\ \,{\rm for}\, P \in C\, ,$$
$$dg^{*}(P_{12}) = \left( \begin{array}{cc} \delta^{-1} & 0\\ 
0 & \delta^{2} \end{array}\right)\ \, ,\, dg^{*}(P_1) = \left( \begin{array}{cc} \delta^{-2} & 0\\ 
0 & \delta^{3} \end{array}\right)\ \, ,\, dg^{*}(P_0) = \left( \begin{array}{cc} \delta^{-1} & 0\\ 
0 & \delta^{2} \end{array}\right)\ \, ,\, $$
$$dg^{*}(P_{45}) = \left( \begin{array}{cc} \delta^{-1} & 0\\ 
0 & \delta^{2} \end{array}\right)\ \, ,\, dg^{*}(P_{56}) = \left( \begin{array}{cc} \delta^{-2} & 0\\ 
0 & \delta^{3} \end{array}\right)\ \, ,$$ 
$$dg^{*}(P_{67}) = \left( \begin{array}{cc} \delta^{-3} & 0\\ 
0 & \delta^{4} \end{array}\right)\ \, ,\, 
dg^{*}(P_{7}) = \left( \begin{array}{cc} \delta^{-4} & 0\\ 
0 & \delta^{5} \end{array}\right)\ \, .$$
In particular, at any point $P \in \mathcal F_C$, 
the eigenvalues of $dg^{*}(P)$ are not multiplicatively 
independent, so that $g$ has no Siegel disk at $P \in \mathcal F_C$.
\end{lemma}
\begin{proof} As $g^{*} \vert NS(S) = id_{NS(S)}$, it follows that $g(C_k) = C_k$ for each $k$. In particular $g(P_{ij}) = P_{ij}$. 
Thus $g\vert C_3 = id_{C_{3}}$, as $C_3 \simeq \bP^1$ 
and $g \vert C_3$ fixes the three points 
$C_3 \cap C_2$, $C_3 \cap C_0$, $C_3 \cap C_4$ on $C_3$. Then, 
$dg^{*}(P)$ is as claimed for $P \in C_3$, by $g^{*}\sigma_S = \delta\sigma_S$. In particular, $d(g \vert C_2)^{*}(P_{23})$ is the multiplication by $\delta$. 
Hence $d(g \vert C_2)^{*}(P_{12})$ is the multiplication by $\delta^{-1}$ (and $g\vert C_2$ has no other fixed point), as $C_2 \simeq \bP^1$. Thus, by $g^{*}\sigma_S = \delta \sigma_S$ again, it follows that 
$d(g \vert C_1)^{*}(P_{12})$ is the multiplication by $\delta^2$. Then, as $C_1 \simeq \bP^1$, the automorphism $g \vert C_1$ has one more fixed point, say, $P_1$, and $d(g \vert C_1)^{*}(P_1)$ is the multiplication by $\delta^{-2}$. 
Then again by $g^*\sigma_S = \delta\sigma_S$, one can diagonalize the action 
$dg^{*}(P_1)$ as claimed. In this way, we figure up the set ${\mathcal F}_{C}$ and the induced actions on the cotangent spaces as claimed. Note 
that, by definition, $\delta^{a}$ and $\delta^{b}$ ($a, b \in \bZ$) are not multiplicatively independent. From this, the last statement follows. 
\end{proof}
Let us define the rational function $\gamma(x) \in \bQ(x)$ by
$$\gamma(x) = \frac{1 + x - x^3 - x^4 - x^5 - x^6 - x^7 - x^8 + x^{10} + x^{11}}{1 + x - x^3 - x^4 - x^5 - x^6 - x^7 + x^{9} + x^{10}}\, .$$
We note that the denominator and the numerator are reciprocal of 
degree $10$ and of degree $11$. Thus, 
$\gamma(x)$ is also written in the form
$$\gamma(x) = \frac{(x+1)f_{1}(x + \frac{1}{x})}{f_{2}(x + \frac{1}{x})}\, ,$$
where $f_1(t)$ and $f_2(t)$ are some polynomials of degree $5$ 
with rational coefficients.
\begin{lemma}\label{fixpt} $g$ has one more fixed point 
$Q \in S \setminus \cup_{k=0}^{7} C_k$. Moreover, the action $dg^{*}(Q)$ is diagonalized as follows:
$$dg^{*}(Q) = \left( \begin{array}{cc} \epsilon_1 & 0\\ 
0 & \epsilon_2 \end{array}\right)\ \, .$$
Here $\epsilon_1$ and $\epsilon_2$ are the roots of the quadratic equation
$$x^2 - \gamma(\delta)x + \delta = 0\, .$$
Their approximate values are 
$$\epsilon_1 = -(0.8886...) -(0.45858...)i\, ,\, \epsilon_2 = -(0.94351...) - 
(0.33133...)i\, .$$ 
\end{lemma}
\begin{proof} As $S$ contains no irreducible complete curve other than $\{C_k\}_{k=0}^{7}$, the fixed points outside $\cup_{k=0}^{7}C_k$ are all isolated 
and finite (even if they exist). 
Let $t \ge 0$ be the number of the fixed points off $\cup_{k=0}^{7} C_k$, 
counted with multiplicities. 
\par
\vskip 4pt
Let us determine $t$ first. By the topological Lefschetz fixed point formula (see eg. \cite{GMa93}, Theorem (10.3)), we have that 
$$\sum_{F} n(F) = 
\sum_{k=0}^{4} (-1)^{k} {\rm tr}\,(g^{*} \vert H^{k}(S, \bZ))\, .
$$
Here the sum in the left hand side runs over all the irreducible components 
of $S^g$. For an isolated point $P$, $n(P)$ is the multiplicity and 
for a fixed smooth curve $C$, the number $n(C)$ is the topological Euler number of $C$ if it is smooth, of multiplicity $1$ (see ibid.). In our case, each irreducible component of $\mathcal F_C$ is smooth, of multiplicity $1$ by the explicit description in Lemma (\ref{fixcurve}). Thus,
$$\sum_{F} n(F) = 7 + 2 + t = 9 +t\, .$$
Here $2$ is the topological Euler number of the fixed curve 
$C_3 \simeq \bP^1$. 
On the other hand, using the fact that $S$ is a K3 surface and the fact that 
the characteristic polynomial of $g^{*} \vert H^2(S, \bZ)$ is $(x-1)^8\varphi_{14}(x)$, we can calculate the right hand side as follows:
$$\sum_{k=0}^{4} (-1)^{k} {\rm tr}\,(g^{*} \vert H^{k}(S, \bZ)) 
= 1 + 1 + {\rm tr}\, ((x-1)^8\varphi_{14}(x))$$
$$= 1 + 1 + {\rm tr}\, ((x-1)^8) + {\rm tr}\, (\varphi_{14}(x)) = 1 + 1 + 8 + 
0 = 10\, .$$
Thus, $9 + t = 10$ and $t = 1$. 
\par
\vskip 4pt 
Hence the fixed point outside 
$\cup_{k=0}^{7} C_k$ is just one point with multiplicity $1$. 
We denote this point by $Q$. Let $\epsilon_1$ and $\epsilon_2$ be the eigenvalues of $dg^{*}(Q)$. As $Q$ is an isolated fixed point 
of multiplicity $1$, we have 
$$\epsilon_1 \not= 1\, ,\, \epsilon_2 \not= 1\, .$$
\par
\vskip 4pt 
Let us determine $\epsilon_1$ and $\epsilon_2$. First of all, 
by $g^{*}\sigma_S = 
\delta \sigma_S$, we have
$$\epsilon_1 \epsilon_2 = \delta\, .$$
Next let us compute the sum $\epsilon_1 + \epsilon_2$. 
For this aim, we want to apply an appropriate form of holomorphic 
Lefschetz fixed point formula. 
In our case, 
$g$ has a fixed curve and $g$ is of infinite order. So, we {\it can not 
directly apply} Atiyah-Bott's one \cite{AB68} or Atiyah-Singer's one 
\cite{AS68}.
On the other hand, $S^g$ is smooth and of multiplicity $1$ at each irreducible component. Thus, we {\it can apply} Toledo-Tong's form of the holomorphic Lefschetz fixed point formula (\cite{TT78}, the formula (*) in Page 519 or Theorem (4.10), applied for $E = \mathcal O_S$):
$$\sum_{F} L(F) =  \sum_{k=0}^{2} (-1)^k{\rm tr}\, 
(g^{*}\vert H^k(\mathcal O_S))\, .$$
Here the sum in the left hand side runs over all the irreducible components 
of $S^g$. The local contribution terms $L(F)$ are calculated 
as follows (See ibid.). 
For an isolated point $P$, 
$$L(P) = \frac{1}{(1-\alpha)(1-\beta)}\, .$$
Here $\alpha$ and $\beta$ are the eigenvalues of $dg^{*}(P)$. 
For a smooth curve $F$ for which $g^{*}\vert N^{*} = \lambda$ 
on the conormal bundle $N^{*}$, we have
$$L(F) = \int_{F} {\rm Td}(F)\cdot \{{\rm ch}(\mathcal O_F) - 
\lambda \cdot {\rm ch}(N^{*})\}^{-1}$$
$$= \int_{F} (1 - \frac{c_1(F)}{2})\cdot \frac{1}{1-\lambda} \cdot 
(1 + \frac{\lambda}{1-\lambda}N^{*}) = \frac{1-p_a(F)}{(1 -\lambda)} + 
\frac{\lambda \cdot {\rm deg} N^{*}}{(1-\lambda)^2}\, .$$
For our $C_3 \simeq \bP^1$, it is
$$L(C_3) = \frac{1}{(1-\delta)} + \frac{\delta \cdot 2}{(1-\delta)^2} = \frac{1 + \delta}{(1- \delta)^2}\, .$$
Thus, the left hand side for our $(S, g)$ is:
$$\sum_{F} L(F) = \frac{3}{(1 - \delta^{-1})(1- \delta^2)} 
+ \frac{2}{(1 - \delta^{-2})(1- \delta^3)} + 
\frac{1}{(1 - \delta^{-3})(1- \delta^4)}$$ 
$$+ \frac{1}{(1 - \delta^{-4})(1- \delta^5)} 
+\frac{1+\delta}{(1 - \delta)^2} 
+ \frac{1}{(1 - \epsilon_1)(1- \epsilon_2)}\, .$$
Let us compute the right hand side. 
$g^{*} \vert H^2(\mathcal O_S)$ is the multiplication by 
$\delta^{-1}$. This is because $H^2(\mathcal O_S)$ is the Serre dual of 
$H^0(\Omega_S^2) = \bC \sigma_S$ and 
$g^{*}\sigma_S = \delta\sigma_S$. 
Thus 
$$\sum_{k=0}^{2} (-1)^k{\rm tr}\, 
(g^{*}\vert H^k(\mathcal O_S)) = 1 + \frac{1}{\delta}\, .$$
Hence
$$1 + \frac{1}{\delta} = \frac{3}{(1 - \delta^{-1})(1- \delta^2)} 
+ \frac{2}{(1 - \delta^{-2})(1- \delta^3)} 
+ \frac{1}{(1 - \delta^{-3})(1- \delta^4)}$$ 
$$+ \frac{1}{(1 - \delta^{-4})(1- \delta^5)} 
+\frac{1+\delta}{(1 - \delta)^2} + \frac{1}{(1 - \epsilon_1)(1- \epsilon_2)}\, .$$
By transposition, we can rewrite this equation in the following form:
$$\frac{1}{(1 - \epsilon_1)(1 -\epsilon_2)} = f(\delta)\, .$$
To get an explicit form of $f(\delta)$, {\it we regard as if $\delta$ is an indeterminate element} and use Mathematica program, Together. The result is:
$$f(\delta) = \frac{1 + \delta - \delta^3 - \delta^4 - \delta^5 
- \delta^6 - \delta^7 + \delta^9 + \delta^{10}}
{(-1 + \delta)^2\delta(1+\delta)(1+\delta + \delta^2)(1+\delta + \delta^2 
+ \delta^3 + \delta^4)}\, .$$
Note that 
$$(1-\epsilon_1)(1-\epsilon_2) = 1 + \delta - (\epsilon_1 + \epsilon_2)\, .$$ 
Then, by the formula above, we obtain that 
$$\epsilon_1 + \epsilon_2 = 1 + \delta - \frac{1}{f(\delta)}\, .$$
In order to simplify the right hand side, we again regard $\delta$ as an indeterminate element and use Mathematica program, Together. The result is: 
$$1 + \delta - \frac{1}{f(\delta)} = \gamma(\delta)\, .$$
Here $\gamma(x)$ is the rational function defined just before Lemma 
(\ref{fixpt}).
Thus, 
$$\epsilon_1 + \epsilon_2 = \gamma(\delta)\, .$$ 
Hence $\epsilon_1$ and $\epsilon_2$ are the roots of the quadratic equation
$$x^2 - \gamma(\delta)x + \delta = 0\, ,$$
as claimed. Now using Mathematica again, we can find approximate values of $\epsilon_1$ 
and $\epsilon_2$ as follows. First, by substituting 
$$\delta = (-0.9903988352300419...) - (0.13823945592693967...)i$$ 
into $\gamma(\delta)$ by Mathematica, we obtain
$$\gamma(\delta) = (0.0548626217729844...) - (0.7899228027367716...)i\,.$$
Our quadratic equation is then
$$x^2 -  ((0.0548626217729844...) - (0.7899228027367716...)i)x$$
$$ + 
((-0.9903988352300419...) - (0.13823945592693967...)i) = 0\, .$$
Solving this equation by Mathematica program, NSolve, we obtain 
approximate values of 
$\epsilon_1$ and $\epsilon_2$ as claimed. They are certainly different and therefore $dg^{*}(Q)$ can be diagonalized. 
\end{proof}
\begin{lemma}\label{galois} $g$ has a Siegel disk at $Q$.
\end{lemma}
Note that the diagonalization in Lemma (\ref{fixpt}) is just on 
the cotangent space level and far from local coordinate level. 
\begin{proof} It suffices to check that 
$\epsilon_1$ and $\epsilon_2$ in Lemma (\ref{fixpt}) satisfy the conditions 
(1)-(3) in Theorem (\ref{siegeldisk}). 
\par
\vskip 4pt
(1) is clear as $\epsilon_1$ and $\epsilon_2$ are the roots of $x^2 - \gamma(\delta) x + \delta = 0$, and both $\delta$ and $\gamma(\delta)$ are algebraic numbers. 
\par
\vskip 4pt
Let us check (2). Mathematica program, Abs applied for $\epsilon_1$ and 
$\epsilon_2$ certainly indicates the result. However, to conclude that some value $x$ is {\it exactly} $1$, computation based on approximate values of 
$x$ seems 
insufficient. Here is a safer argument. Consider 
$$e_1 := \frac{\epsilon_1^2}{\delta}\, ,\, 
e_2 := \frac{\epsilon_2^2}{\delta}\, .$$
As $\vert \delta \vert = 1$, it suffices to show that $\vert e_1 \vert = \vert e_2 \vert = 1$.  We have $e_1e_2 = 1$ by $\epsilon_1\epsilon_2 = \delta$. 
We also have 
$$e_1 + e_2 = 
\frac{(\epsilon_1 + \epsilon_2)^2 - 2\epsilon_1\epsilon_2}{\delta} 
= \frac{\gamma(\delta)^2}{\delta} -2\, .$$
Recall the second expression of $\gamma(\delta)$ given 
just before Lemma (\ref{fixpt}): 
$$\gamma(\delta) = \frac{(\delta +1)f_1(\delta + \frac{1}{\delta})}
{f_2(\delta + \frac{1}{\delta})}\, .$$
Then, $e_1 + e_2 = k(\delta)$, where
$$k(\delta) = \frac{(\gamma(\delta))^2}{\delta} -2 = 
(\delta + \frac{1}{\delta} + 2)
\{\frac{f_1(\delta + \frac{1}{\delta})}{f_2(\delta + \frac{1}{\delta})}\}^2 -2
\, ,$$
and $e_1$ and $e_2$ are the roots of the quadratic equation
$$x^2 - k(\delta)x + 1 = 0\, .$$
By the quadratic formula, we have  
$$e_1, e_2 = \frac{k(\delta) \pm \sqrt{k(\delta)^2 -4}}{2}\, .$$ 
Here, $k(\delta)$ is {\it real}, as 
$$\delta + \frac{1}{\delta} = \delta + \overline{\delta}$$
is real by $\vert \delta \vert = 1$. Thus, $\vert e_1 \vert = \vert e_2 \vert = 1$ if and only if $\vert k(\delta) \vert \le 2$ from the quadratic formula 
above. Substituting 
$$\delta = -(0.9903...) - (0.3182...)i$$ 
into $k(\delta)$ by Mathematica, 
we have
$$k(\delta) = \frac{\gamma(\delta)^2}{\delta} - 2 
 = -(1.3730...) + (9.4799... \times 10^{-17})i\, .$$
Thus $\vert k(\delta) \vert < 2$. 
\par
\vskip 4pt
We should remark that there appears an error term 
$$(9.4799... \times 10^{-17})i\, .$$ 
in the above expression of $k(\delta)$. However, {\it this does not matter}, 
because it is extremely small compared with the real part 
and we {\it know} that $k(\delta)$ is certainly real. 
Hence the assertion 
(2) holds for our $\epsilon_1$ and $\epsilon_2$. 
\par
\vskip 4pt
It remains to check (3) for our $\epsilon_1$ and $\epsilon_2$. Suppose 
$$\epsilon_1^m\epsilon_2^n = 1$$
for $(m, n) \in \bZ^2$. Note that $\delta = \beta_1$ and $\theta = \beta_2$ 
are Galois conjugate, as both are roots of $\varphi_{14}(x) = 0$ (and 
$\varphi_{14}(x)$ is irreducible). Thus, by taking Galois conjugate, we have
$$(\epsilon_1')^m(\epsilon_2')^n = 1\, .$$
Here $\epsilon_1'$ and $\epsilon_2'$ are the roots of the quadratic equation
$$x^2 - \gamma(\theta)x + \theta = 0\, .$$
As 
$$\theta = -(0.371932997164175...) - (0.92825957879273...)i\, ,$$
we have 
$$\gamma(\theta) = (1.495690836752066...) - (2.210575209107991...)i\, ,$$
by Mathematica. Substituting these values into the quadratic equation above and using Mathematica program, NSolve, we find (up to order) that
$$\epsilon_1' = (0.25262...) - (0.37337...)i\,\, ,\,\, 
\epsilon_2' = (1.2430...) -(1.837...)i\, .$$
Clearly $\vert \epsilon_2' \vert > 1$. Thus, from 
$$1 = \vert (\epsilon_1')^m(\epsilon_2')^n \vert = 
\vert (\epsilon_1'\epsilon_2')^m
(\epsilon_2')^{n-m} \vert = \vert (\theta)^{m}(\epsilon_2')^{n-m}\vert = \vert(\epsilon_2')
\vert^{n-m}\, ,$$
we conclude $n = m$. 
Substituting this into $\epsilon_1^n \epsilon_2^m= 1$, 
we obtain
$$1 = (\epsilon_1\epsilon_2)^n = \delta^n\, .$$
Here $\delta$ is not root of unity. This is because the Salem 
number $\alpha_{14} >1$ is a Galois conjugate of 
$\delta$. Hence $n = 0$, and therefore, $m = n = 0$ by $n=m$. This shows (3).
\end{proof}
Now the following Lemma completes the proof of Theorem (\ref{main}):
\begin{lemma}\label{autom} 
${\rm Aut}\, S = \langle g \rangle \simeq \bZ$.
\end{lemma}
\begin{proof} As $\delta$ is not a root of unity, $g$ is of infinite order. 
So, it suffices to show that ${\rm Aut}\, S$ is generated by $g$. 
Let $f \in {\rm Aut}\, S$. As the dual graph of the curves $\{C_k\}_{k=0}^{7}$ 
is the Dynkin diagram $E_8(-1)$ and it has no symmetry, 
we have $f(C_k) = C_k$ ($0 \le k \le 7$). Hence $f^{*} \vert NS(S) = 
id_{NS(S)}$, 
as $\{C_k\}_{k=0}^{7}$ generates $NS(S)$. The natural representation 
of ${\rm Aut}\, S$ on $T(S)$
$$r_T : {\rm Aut}\,S \longrightarrow {\rm O}(T_S)$$
is then injective, as so is on ${\rm O}\,(H^2(S, \bZ))$ (see eg. \cite{BHPV04}, Page 333, Corollary (11.4)). Moreover, as 
$NS(S) \simeq E_8(-1)$ is negative definite, ${\rm Im}\, r_T$ 
is isomorphic to $\bZ$. This is a special case of \cite{Og08}, Theorem (1.5). 
Hence, ${\rm Aut}\, S$ is isomorphic to $\bZ$ as well. 
Let $h$ be a generator of ${\rm Aut}\, S$. 
By replacing 
$h$ by $h^{-1}$ if necessary, we can write $g = h^n$ for some positive integer 
$n$. Let $\varphi(x)$ be the characteristic polynomial of $h^{*} \vert T(S)$. 
As $NS(S)$ is negative definite, $\varphi(x)$ is again a Salem polynomial of degree $14$ (\cite{Og08}, Theorem (3.4)). Let $\beta_{14}$ be the Salem number of $\varphi(x)$. Then, by $g = h^n$, we have
$$\alpha_{14} = \beta_{14}^n\, .$$
On the other hand, $\alpha_{14}$ is the {\it smallest} Salem number of
 degree $14$, as explained in Section 2. Hence 
$n =1$, i.e., $g = h$. 
\end{proof} 

\begin{remark}\label{alternate} Let us consider the pair $(S', g')$ in Remark 
(\ref{K3'}). Then, as Lemmas (\ref{fixcurve}), 
(\ref{fixpt}), we have a similar description of the fixed point set: 
$$(S')^{g'} = C_3' \cup \{P_1', P_{12}', P_0', P_{45}', P_{56}', 
P_{67}', P_7' \} \cap \{Q'\}\, .$$
{\it However, $g'$ has no Siegel disk.} In fact,  
The eigenvalues of $d(g')^{*}(P')$ ($P' \in C_3'$), $d(g')^{*}(P_i')$, 
$d(g')^{*}(P_{ij}')$ are the same as the eigenvalues of $dg^{*}(P)$ ($P \in C_3$), $dg^{*}(P_i)$ and $dg^{*}(P_{ij})$, and they are not multiplicatively independent. 
The eigenvalues 
of 
$d(g')^{*}(Q')$ are $\epsilon_1'$ and $\epsilon_2'$. Here, $\epsilon_1'$ and 
$\epsilon_2'$ are 
the numbers defined at the last part of the proof of Lemma (\ref{galois}). 
Then $\vert \epsilon_1' \vert > 1$ as observed there. So, 
$g'$ has no Siegel disk at $Q'$, either.
\end{remark}
\begin{remark}\label{4thsalem} Recall (from Section 2 (i)) that 
$$\Phi_{14}(x) = x^{14} -x^{12} -x^{7} -x^{2} +1$$
is the Salem polynomial of the $4$-th smallest known Salem number 
$$A_{14} = 1.20261...\,\, .$$
$\vert \Phi_{14}(\pm 1) \vert = 1$, and $\Phi_{14}(x) = 0$ 
has two particular roots on the unit circle
$$\delta' := -(0.45829...) - (0.88799...)i\, ,$$
$$\theta' := -(0.96815...) - (0.25034...)i\, .$$
Then, starting from $\delta'$ and arguing exactly in the same way as 
in Theorem (\ref{main}), we also obtain a K3 surface automorphism of 
topological entropy $\log\, A_{14}$, with one pointwisely fixed 
smooth rational curve and a Siegel disk. The resultant action on 
the Siegel disk is given by ${\rm diag}\, (\rho_1, \rho_2)$, where
$$\rho_1 = -(0.29457...) - (0.95562...)i\, ,\, 
\rho_2 = (0.98436...) - (0.17614...)i\, .$$
\end{remark}
\section{Proof of Theorem (\ref{main2})}
\setcounter{lemma}{0}
\noindent
In this section, we shall prove Theorem (\ref{main2}). Let $S$ be an Enriques surface and $g$ be an autmorphism of $S$. Let us denote the free part of $H^2(S, \bZ)$ by $L$. Then, by \cite{BHPV04}, Page 339, 
Lemma (15.1) (iii), we have
$$L \simeq H \oplus E_{8}(-1) \simeq \, E_{10}(-1)\, .$$
Here the last isomorphism comes from the fact that both $H \oplus E_{8}(-1)$ 
and $E_{10}(-1)$ are even unimodular lattices of signature $(1,9)$. In fact, we can then apply \cite{Se73}, Page 54, Theorem 5.
We denote by $L(2)$ the lattice such that $L(2) = L$ as $\bZ$-module and 
$$(x, y)_{L(2)} := 2(x, y)_{L}$$ 
for each $x, y \in L(2) = L$. Note that $g^{*}$ 
is also an automorphism of the new lattice $L(2)$. Then $g^{*}$ acts on the discriminant group 
$$A_{L(2)} := L(2)^{*}/L(2) \simeq L/2L \simeq \bF_2^{10}\, .$$
Let $\overline{\Phi_{g^{*}}}(x)$ be the characteristic polynomial of $g^{*} \vert A_{L(2)}$. Then $\overline{\Phi_{g^{*}}}(x) \in \bF_2[x]$. More precisely, 
$\overline{\Phi_{g^{*}}}(x)$ is the mod $2$ reduction of the characteristic polynomial $\Phi_{g^{*}}(x)$ of $g^{*} \vert L = g^{*} \vert L(2)$. 
\begin{lemma}\label{order} No irreducible component of 
$\overline{\Phi_{g^{*}}}(x) \in \bF_2[x]$ is of degree $5$.
\end{lemma} 
\begin{proof}
Let $\pi : \tilde{S} \longrightarrow S$ be the universal cover of $S$. 
Then $\tilde{S}$ is a K3 surface and $\pi$ is of degree $2$ 
(See eg. \cite{BHPV04}, Page 339, Lemma (15.1)(ii)). 
We denote by $\iota$ the covering involution of $\pi$. Following \cite{Nm85}, Page 203, line 6, we define
$$M := \{x \in H^2(\tilde{S}, \bZ)\, \vert \iota^*x = x\}\, ,\, N := \{x \in H^2(\tilde{S}, \bZ)\, \vert \iota^*x = -x\}\, .$$
By \cite{Nm85}, Proposition (2.3), we have
$$M = \pi^{*}(L) = L(2)\, .$$
Here the last equality is nothing but the following obvious relation 
$$(\pi^*x, \pi^*y)_{\tilde{S}} = 2(x, y)_{S}\,\, \forall x, y \in L\,.$$
Let $A_{M} = M^{*}/M$ and $A_{N} = N^{*}/N$ be the discriminant groups 
of $M$ and $N$. Note that  
$$A_N \simeq A_M = A_{L(2)} \simeq \bF_2^{10}\, .$$
Here the first isomorphism is given by the natural surjective morphisms: 
$$H^2(\tilde{S}, \bZ) \longrightarrow M^{*}/M\, ;\, x \mapsto (x, *)\, {\rm mod}\, M\, ,$$
$$H^2(\tilde{S}, \bZ) \longrightarrow N^{*}/N\, ;\, x \mapsto (x, *)\, {\rm mod}\, N .$$
We note that these two morphisms are certainly surjective as $H^{2}(\tilde{S}, \bZ)$ is unimodular and both $N$ and $M$ are primitive. Then, it is an easy fact that the both kernels are $N \oplus M$. This is a special case of 
\cite{Ni80}, Corollary (1.5.2). 
\par
\vskip 4pt
Let $\tilde{g} \in {\rm Aut}\, \tilde{S}$ be one of the two possible lifts of 
$g$ on $\tilde{S}$. Then $g \circ \pi = \pi \circ \tilde{g}$ and 
$\tilde{g} \circ \iota = \iota \circ \tilde{g}$. Thus, $\tilde{g}^{*}$ 
preserves both $M$ and $N$. Hence, $\tilde{g}^{*}$ induces actions on 
$M$ and $N$, 
and consequently, on the discriminant groups $A_{M}$ and $A_{N}$. 
Moreover, under the isomorphism of discriminant groups above, we have
$$\tilde{g}^{*} \vert A_{N} = \tilde{g}^{*} \vert A_{M} = 
g^{*} \vert A_{L(2)}\, .$$  
Here the last equality follows from $\pi \circ \tilde{g} = g \circ \pi$. 
Thus, the characteristic polynomial of $g^{*} \vert A_{N}$ is the same as the characteristic polynomial $\overline{\Phi_{g^{*}}}(x)$ of 
$g^{*} \vert A_{L(2)}$. So, if 
$\overline{\Phi_{g^{*}}}(x)$ would have an irreducible factor of degree $5$, then the corresponding eigenvalue of $\tilde{g}^{*} \vert A_N$ 
would be an element of 
$\bF_{32} \setminus \bF_2$. Here $32 = 2^{5}$. As 
$$(\bF_{32})^{\times} \simeq \bZ/(32-1)\bZ = \bZ/31\bZ\, ,$$
the order of $\tilde{g}^{*} \vert A_{N}$ would then 
be divisible by $31$. Thus, the order 
of $\tilde{g}^{*} \vert N$, which is actually finite 
(Lemma (\ref{trans}) below), 
would be also divisible by $31$. However, this is impossible by 
the next a bit more precise Lemma (\ref{trans}).  
\end{proof}
\begin{lemma}\label{trans} Under the same notations as in the proof of Lemma 
(\ref{order}), the order of $\tilde{g}^{*} \vert N$ is finite, say, $d$. 
Let 
$$d = \Pi_{k=1}^{n} p_{k}^{m_k}$$
 be the prime decomposition of $d$. Then 
each primary factor $p_{k}^{m_k}$ belongs to
$$\{2, 2^2, 2^3, 2^4, 3, 3^2, 5, 7, 11, 13\}\,.$$ 
\end{lemma}

\begin{proof} As the lattice 
$M = \pi^{*}L$ is of signature $(1, 9)$ with pure Hodge type $(1,1)$, the lattice $N$ is of signature $(2, 10)$ and $N$ admits the following real 
Hodge decomposition:
$$N_{\bR} = Q \oplus P\, .$$
Here  
$$Q := N_{\bR} \cap H^{1,1}(\tilde{S})\,\, ,\,\, P := \bR \langle {\rm Re}\, 
\sigma_{\tilde{S}}, {\rm Im}\, \sigma_{\tilde{S}} \rangle,$$ 
and $\sigma_{\tilde{S}}$ is a nowhere vanishing global holomorphic 
$2$-form on $\tilde{S}$. 
As $\tilde{g}^{*} \vert N$ preserves the Hodge decomposition, we have 
$$\tilde{g}^{*} \vert N \in {\rm O}(P) \times {\rm O}(Q)\, .$$
Here $P$ is positive definite and $Q$ is negative definite. Hence $\tilde{g}^{*} \vert N$ is diagonalizable and the eigenvalues are of absolute value $1$. 
On the other hand, $\tilde{g}^{*} \vert N$ is defined over $\bZ$. Thus, 
all the eigenvalues are roots of unity (Kronecker's theorem). 
Hence $\tilde{g}^{*} \vert N$ is of finite order, say, $d$. We denote the prime decomposition of $d$ as in the statement. Then $(\tilde{g}^{*} \vert N)^{e}$ with $e = d/p_{k}^{m_k}$, is of order $p_{k}^{m_k}$. As $(\tilde{g}^{*} \vert N)^{e}$ 
is defined over $\bZ$, all primitive $p_k^{m_k}$-th roots of unity appear 
as eigenvalues of $(\tilde{g}^{*} \vert N)^{e}$. As easily seen, their cardinality is exactly $p_k^{m_k -1}(p_k -1)$. As ${\rm rank}\, N = 12$, this number can not exceed $12$, that is,  
$$p_k^{m_k -1}(p_k -1) \le 12\, .$$
Solving this inequality, we obtain the result.
\end{proof}  

Now we are ready to complete the proof of Theorem (\ref{main2}). If $h(g)$ 
would be the logarithm of the Lehmer number, then, as $L$ is of rank $10$, the characteristic polynomial $\Phi_{g^{*}}(x)$ of $g^{*} \vert L$ 
would be the Lehmer polynomial: 
$$\varphi_{10}(x) = x^{10} +x^{9} -x^{7} -x^{6} -x^{5} -x^{4} -x^{3} + x + 1\, .$$
However, this is impossible by Lemma (\ref{order}) and the following:
\begin{lemma}\label{mod2} Let $\overline{\varphi_{10}}(x) \in \bF_2[x]$ 
be the mod $2$ reduction of the Lehmer polynomial. 
Then the irreducible decomposition of $\overline{\varphi_{10}}(x)$ is:
$$\overline{\varphi}_{10}(x) = 
(x^5 + x^3 +x^2 + x+1)(x^5 + x^4 +x^3 + x^2 +1)\, .$$
\end{lemma}
\begin{proof}
This immediately follows from Mathematica program, Factor:
$${\rm Factor}\, [x^{10} +x^{9} -x^{7} -x^{6} -x^{5} -x^{4} -x^{3} + x + 1, {\rm Modulus} \rightarrow 2]\, .$$
\end{proof}
This completes the proof of Theorem (\ref{main2}). 

\begin{remark}\label{5th}
The second, third, forth smallest known Salem numbers are of degree $> 10$. 
So, they can not be realized as the antilogarithm of the topological entropy of an Enriques surface automorphism. Recall (from Section 2 (i)) that 
$$\Phi_{10}(x) = x^{10} -x^{6} -x^{5} - x^{4} +1\, $$
is the Salem polynomial of the fifth smallest known Salem number
$$A_{10} = 1.216391...\, .$$
Again by Mathematica program, Factor, applied for the mod $2$ reduction 
$\overline{\Phi_{10}}(x)$ of $\Phi_{10}(x)$, we obtain the irreducible 
factorization:
$$\overline{\Phi_{10}}(x) = 
(x^5 + x^4 + x^2 + x + 1)(x^5 + x^4 + x^3 + x +1)\, .$$
Thus, the fifth smallest known Salem number $A_{10}$ can not 
be realized, either. In conclusion, {\it none of the smallest five 
Salem numbers 
(listed in Section 2 (i)) can be realized as the 
antilogarithm of the topological entropy of 
an Enriques surface automorphism.}  
\end{remark}

We conclude this note by the following natural, probably tractable, open problems relevant to our Theorems (\ref{main}), (\ref{main2}) with a few remarks:
\begin{question}\label{open}

(1) What is the smallest Salem number that can be realized as the antilogarithm of the topological entropy of a K3 surface automorphism? 

(2) Is there a surface automorphism having more than one Siegel disks 
at the same time? 

(3) What is the smallest Salem number that can be realized as the antilogarithm of the topological entropy of an Enriques surface automorphism? 

(4) What is the smallest Salem number that can be realized as the antilogarithm of the topological entropy of an automorphism of a {\it generic} Enriques 
surface? 

(5) Are there compact hyperk\"ahler manifolds (of dimension $2n$) having 
automorphisms with $2n$-dimensional Siegel disks? 
\end{question}

\begin{remark}\label{openrm} Here are a few remarks about some of questions 
above.

{\it For Question (1).} As the Salem numbers in question are of degree $\le 22$, we only need to see the realizability of the first and second smallest Salem numbers $\alpha_{10}$ and $\alpha_{18}$. 

{\it For Question (4).} The automorphism 
group of a generic Enriques surface is isomorphic to the 
$2$-congruence subgroup ${\rm O}^{+}(E_{10}(-1))(2)$ of 
${\rm O}^{+}(E_{10}(-1))$. This is proved by \cite{BP83}, Theorem (3.4) 
and \cite{Nm85}, Theorem (5.10). (See also the precise meaning ``generic" there.) Thus, this is also a purely group theoretical problem.

{\it For Question (5).} In the terminology of \cite{Be83}, 
Page 759, Th\'eor\`eme, a compact hyperk\"ahler manifold of dimension 
$2n$ is a Ricci flat compact K\"ahler manifold with ${\rm Sp}(n)$ holonomy and a Calabi-Yau manifold of dimension $m \ge 3$ is a Ricci flat compact K\"ahler manifold with ${\rm SU}(m)$ holonomy. Let $X$ be a Calabi-Yau manifold of dimension $m \ge 3$. Then $X$ is projective (\cite{Be83}, Page 760, Proposition 1), and therefore, the action of 
${\rm Aut}\, X$ on the space of holomorphic $m$-forms is finite cyclic 
(\cite{Ue75}, Page 178, Proposition 14.5). So, Calabi-Yau manifolds of dimension $m \ge 3$ can not admit automorphisms with $m$-dimensional Siegel disks. For the same reason, compact 
hyperk\"ahler manifolds having automorphisms with $2n$-dimensional Siegel disks
can not be projective as well. A bit more precisely, they are in fact of 
algebraic dimension $0$ (\cite{Og08}, Theorem (2.4)). See also \cite{Og09} for the explicit description of the topological entropy of automorphisms of compact hyperk\"ahler 
manifolds and \cite{Zh08} for a more algero-geometric aspect of 
the topological entropy of automorphisms of higher dimensional manifolds.
\end{remark}

%\newpage

\vskip .2cm \noindent
Keiji Oguiso \\
Department of Mathematics, Osaka University\\
Toyonaka 560-0043 Osaka, Japan\\
oguiso@math.sci.osaka-u.ac.jp

\end{document}